\documentclass[10pt]{amsart}
\usepackage{amsmath}
\usepackage{amsfonts}
\usepackage{amssymb}
\usepackage{amscd}
\usepackage{graphics}
\usepackage{graphicx}
\usepackage{xypic}
 \input{xy}
  \xyoption{all}
\usepackage{xy}
\usepackage{clrscode}
\usepackage{hyperref}
\usepackage{algorithmic}
\newtheorem{theorem}{Theorem}[section]
\newtheorem{definition}[theorem]{Definition}

\newtheorem{example}[theorem]{Example}

\newtheorem*{thm}{Theorem}

\algsetup{linenodelimiter=}

\begin{document}

% Title of document, usually lower case except for first word
% and proper nouns.  Avoid unnecessary symbols.
\title{An algorithm for low dimensional group homology} 

% If the title is too long for the running head, use
% the following command to specify a short title:
%\shorttitle{Shorter title}

\author{Joshua Roberts}             % First Author
\begin{abstract}
Given a finitely presented group $G$, Hopf's formula expresses the second integral homology of $G$ in terms of 
 generators and relators. We give an algorithm that exploits Hopf's formula to estimate $H_2(G;k)$,
 with coefficients in a finite field k, and give examples using $G=SL_2$ over  specific rings of integers. These examples are
 related to a conjecture of Quillen.
\end{abstract}

% Leave these items like this, and the journal will fill them in.
%\received{Month Day, Year}   % receive date (for example: October 11, 1999)
%\revised{Month Day, Year}    % date of revision; omit, if no revision;
                             % if multiple revisions, separate by commas
%\published{Month Day, Year}  % publish date
%\submitted{Graham Ellis}  % Name of Journal's Editor, who handled Article 
%\volumeyear{2008} % Volume Year
%\volumenumber{10} % Volume Number 
%\issuenumber{2}   % Issue Number
%\startpage{1}     % PageNumber of first page
%\webaddress{http://intlpress.com/HHA/v10/n2/a?}
% If copyright is retained by author, comment this out:
%\owner{International Press}

\maketitle

% Text of Document.  Use constructs such as \section, \subsection,
% \begin{theorem} ... \end{theorem}, \begin{proof} ... \end{proof}, etc.
% See 

\section{Introduction}
The purpose of this note is to give an algorithm that allows us to 
estimate the second homology group of any finitely presented group.
More precisely, given a finitely presented group $G$ and a 
finite field $k$, the second homology group $H_2(G;k)$ with coefficients
in $k$ is a finite dimensional vector space over $k$, where $G$ acts on $k$ trivially. Our algorithm 
gives an upper bound for the dimension of $H_2(G;k)$ and, in 
particular cases, the algorithm calculates precisely this dimension.
This algorithm is an improvement of existing algorithms to compute $H_2(G)$.
For instance, the algorithms included in the GAP \cite{GAP4} packages 
``cohomolo" \cite{holt-cohomolo:www} and ``HAP" \cite{ellis-hap:www}
are effective on finite groups and finite or $p$-groups respectively. The algorithm
presented here effectively finds a bound for the homology of any finitely presented group.

A motivational problem for low dimensional group homology is the
study of homology for groups of the form $GL_j(A)$ where $GL_j$ is a
finite rank general linear group scheme and $A$ is a ring of
arithmetic interest. An approach to this problem is to consider the
diagonal matrices inside $GL_j$. Let $D_j$ denote the subgroup
formed by these matrices. Then the canonical inclusions $D_j\subset
GL_j$ for $j=0,1,...$ induce homomorphisms on group homology with
$k$-coefficients
\[
\rho^{A,p}_{i,j}\colon H_i(D_j(A);k)\to H_i(GL_j(A);k)
\]
where $k$ is the field of prime order $p$, $i$ is called the
homological dimension and $j$ the rank. In this context, a
celebrated conjecture of Quillen \cite{quillen:71} implies that
$\rho^{A,p}_{i,j}$ is an epimorphism for $A=\mathbb Z[1/p,\zeta_p]$,
$p$ a regular odd prime, $\zeta_p$ a primitive $p^{\text{th}}$ root of unity 
and any values of $i$ and $j$. For a survey
on the current status of this conjecture we cite \cite{anton:08}.

By a spectral sequence argument applied to the group extension
$$
1\to SL_j(A)\to GL_j(A)\to D_1(A)\to 1
$$
given by the determinant map we can reformulate Quillen's conjecture
in terms of $H_i(SL_j(A);k)$. In the particular case $j=2$ this
homology has been studied extensively by using the theory of
buildings. However, based on this theory we can calculate this
homology only for $i$ sufficiently large \cite{borelserre:76}. The problem of calculating
$H_i(SL_2(A);k)$ in low dimensions turns out to be highly
non-trivial even when $i=2$.

Examples in Section \ref{examples} confirm the results in \cite{anton:08}, as well as give a new finding:
\begin{thm}The dimension of $H_2(SL_2( \mathbb Z[1/7, \zeta_7]); \mathbb F_7)$ as a vector space over $\mathbb F_7$ 
 is at most 6.
\end{thm} 
\noindent Our calculations were done with the computational algebra program GAP  and the GAP commands
to carry out these calculations are given in Appendix \ref{appendix}.

%\ack
The author would like to thank Marian Anton for his enlightening discussions and Mark Dickinson for many helpful comments and recommended revisions.

\section{First Homology Group}

We consider a group given by a finite set of generators and a finite
set of relators. If we denote this group by $G$ then there is a
short exact sequence
$$
1\to R\to F\to G\to 1
$$
where $F$ is a finitely-generated free group and $R$ is a normal
subgroup of $F$ such that if $F$ acts on $R$ by conjugation then $R$
is a finitely-generated $F$-module. Here if $F$ and $R$ are two
groups not necessarily commutative then an $F$-module structure on
$R$ is an assignment $r\mapsto r^f$ for $r\in R$ and $f\in F$ such
that
\begin{eqnarray*}
r^1&=&r\\
(r_1r_2)^f&=&r_1^fr_2^f\\
r^{f_1f_2}&=&(r^{f_1})^{f_2}
\end{eqnarray*}
where, if not otherwise stated, all groups are given
multiplicatively. In this context, it is well known that the first
homology of a group is just another name for its abelianization
\cite{brown:82}. In particular, if we denote by $H_1(G)$ this abelian
group then there is a short exact sequence
\[
1\to R[F,F]\to F\to H_1(G)\to 1
\]
where $[F,F]$ denotes the subgroup of $F$ generated by the
commutators in $F$ and the juxtaposition denotes the operation of
taking the subgroup generated by the parts. Letting $F$ act on
$R[F,F]$ by conjugation, we recognize that $R[F,F]$ is a finitely-generated 
$F$-module. Indeed, the commutator formula
$$
[xy,z]=(xy)^{-1}z^{-1}xyz=y^{-1}x^{-1}z^{-1}xzyy^{-1}z^{-1}yz=[x,z]^y[y,z]
$$
proves that since $F$ is a finitely-generated group then $[F,F]$ is
a finitely-generated $F$-module under conjugation and the same is
assumed about $R$. This argument leads to a deterministic algorithm that gives the structure of
$H_1(G)$. The input is a finite list of generators for $F$, say
$S(F)$, and a finite list of generators for the $F$-module $R$, say
$S(R)$. The output is a list of integers describing the
structure of the finitely-generated abelian group $H_1(G)$.

\subsection{The First Homology Algorithm}

\noindent \textbf{Note on Terminology:} In this paper, we use capital letters, $A$, to denote groups, and $S(A)$ to denote 
a finite set of elements of $A$. However, for simplicity we will simply use capital letters to denote 
these sets in the algorithms.\\

\noindent \textbf{Algorithm 1:} $\proc{FirstHomology}(F,R)$
\begin{algorithmic}[1]
\REQUIRE Free Group $F$, Relators $R$
\ENSURE List of abelian invariants of the finitely presented group $F/R$
\STATE $M:=$ RelationMatrix$(F,R)$
\STATE $N:=$ SmithNormal$(M)$
\RETURN Diagonal$(N)$
\end{algorithmic}
~\\

The GAP command \verb+AbelianInvariants()+ carries out (roughly) the above algorithm. An algorithm
for reducing a matrix to a Smith Normal form is given in \cite{holt:05}. 
Recall that given a finite presentation for $F/R$ that consists of $n$ generators $S(F)$ and $m$ relators $S(R)$, 
there is the associated $n \times m$ relation matrix $M$ whose $(i,j)$ entry is the sum 
of the exponents of all occurrences of the $j$th generator in the $i$th relator.
The resulting list ``Diagonal$(N)$''
is the set of entries in the $(i,i)$ position for $i=1\dots\text{min}(n,m)$ and
consists of positive integers and zeros. The number of zeros is the rank of $H_1(G)$ and 
each positive integer $n$ corresponds to a copy of $\mathbb Z_n$ in the torsion part of $H_1(G)$.

This result can be extended to the case when the homology of $G$ is
taken with trivial coefficients in a finite field say $k$. In this
case, the first homology group of $G$ is denoted by $H_1(G;k)$ and
is a finite dimensional vector space over $k$. Its dimension can be
determined from the universal coefficients \cite{brown:82} short exact
sequence
\[
1\to k\otimes H_1(G)\to H_1(G;k)\to \text{Tor}(H_0(G),k)\to 1
\]
where $H_0(G)$ is the free cyclic group and Tor$(-,k)$ is a functor
vanishing on free abelian groups. The algorithm takes as input the
finite lists $S(F)$ and $S(R)$ from the previous algorithm together
with the characteristic $p$ of the finite field $k$. The output is
an integer representing the dimension of the vector space
$H_1(G;k)$. \\

\subsection{The First Homology with Coefficients Algorithm}
~\\

\noindent \textbf{Algorithm 2:} $\proc{FirstHomologyCoefficients}(F, R, p)$
\begin{algorithmic}[1]
\REQUIRE Free Group $F$, Relators $R$, Prime $p=\text{char}(k)$
\ENSURE Dimension of the vector space $k \otimes H_1(G;k)$ over $k$
\STATE $T:=$FirstHomology$(F,R)$
\STATE $X:=[\ ]$
\FOR{$x \in X$}
\IF{$x \equiv 0 \text{ mod }p$}
\STATE append $x$ to $X$
\ENDIF
\ENDFOR
\RETURN Size($X$)
\end{algorithmic}

\section{Second Homology Group}
Our investigation can be extended to the second homology group of
$G$ which is an abelian group that we denote $H_2(G)$. By a
celebrated formula due to Hopf \cite{brown:82} this group fits into the
following exact sequence
\[
1\to [F,R]\to R\cap [F,F]\to H_2(G)\to 1
\]
where $[F,R]$ is the subgroup of $F$ generated by the commutators
$[f,r]$ with $f\in F$ and $r\in R$. The commutator formula
$$
[x,y^z]=x^{-1}(y^{-1})^zxy^z=x^{-1}z^{-1}y^{-1}zx(yy^{-1})z^{-1}yz=[zx,y][y,z]
$$
proves that $[F,R]$ is a finitely-generated $F$-module under
conjugation. However the intersection $R\cap [F,F]$ is not
determined by any algorithm and we can only estimate the group
$H_2(G)$ as a subgroup of the factor group $R/[F,R]$. This factor
group is abelian since $[F,R]$ contains $[R,R]$ and if we let $F$
act on it by conjugation, this action is trivial. In particular,
since $R$ is a finitely-generated $F$-module it follows that the
factor group $R/[F,R]$ is a finitely-generated abelian group.
Consequently, $H_2(G)$ is a finitely-generated abelian group whose
structure we would like to determine. 

We start with the following exact sequence
\[\label{4}
1\to H_2(G)\to \frac{R}{[F,R]}\to \frac{F}{[F,F]}\to
\frac{F}{R[F,F]}\to 1
\]
in which the last two terms are deterministically determined as
explained above. Moreover, starting with a finite list of generators
$S(R)$ for the $F$-module $R$, we can design a deterministic
algorithm to find a set of generators for $H_2(G)$.

To simplify the discussion, let $k$ denote the finite field of prime
order $p$ and start our investigation with the homology with trivial
coefficients in $k$. By the universal coefficients theorem we have a
short exact sequence
\[
1\to k\otimes H_2(G)\to H_2(G;k)\to \text{Tor}(H_1(G),k)\to 1
\]
whose last term can be determined as follows. For input we start
with the abelian invariants of $H_1(G)$ found by the first algorithm
together with the order $p$ of the field $k$. The output is an
integer say $a$ representing the dimension of the vector space
Tor$(H_1(G),k)$ over $k$. The algorithm is deterministic. \\

\subsection{The Tor Algorithm}
~\\

\noindent \textbf{Algorithm 3:} $\proc{Tor}(F, R, p)$
\begin{algorithmic}[1]
\REQUIRE Free Group $F$, Relators $R$, Prime $p=\text{char}(k)$
\ENSURE Dimension of Tor$(H_1(G),k)$ over $k$
\STATE $A:=$FirstHomology$(F,R)$
\STATE $X:=[\ ]$
\FOR{$x \in A$}
\IF{$x\ne 0$ and $x \equiv 0 \text{ mod }p$ }
\STATE append $x$ to $X$
\ENDIF
\ENDFOR
\RETURN Size$(X)$
\end{algorithmic}
~\\

The first term $k\otimes H_2(G)$ of the exact sequence is a finite
dimensional vector space over $k$ whose dimension can only be
estimated from above by an algorithm that we will
describe next. From exact sequence \eqref{4} we extract the
short exact sequence
\[
1\to H_2(G)\to \frac{R}{[F,R]}\to \frac{R[F,F]}{[F,F]}\to 1
\]
whose last term is a subgroup of the free abelian group $F/[F,F]$.
It is a standard fact that any subgroup of a finitely-generated free abelian group is
free abelian and consequently the above sequence splits. In
particular, by tensoring with $k$ we obtain a short exact sequence
of vector spaces over $k$:
\[
1\to k\otimes H_2(G)\to  k\otimes \frac{R}{[F,R]}\to
 k\otimes \frac{R[F,F]}{[F,F]}\to 1
\]
where the last term can be rewritten as $R[F,F]/R^p[F,F]$. Here
$R^p$ denotes the subgroup of $F$ generated by the $p$-powers of
elements of $R$. In particular, there is a short exact sequence of
finitely-generated abelian groups
\[ \label{8}
1\to
k\otimes\frac{R[F,F]}{[F,F]}\to\frac{F}{R^p[F,F]}\to\frac{F}{R[F,F]}\to
1
\]
whose last two terms are deterministically computable by our first
algorithm. 

\begin{definition}[\cite{faith:00}]
For an abelian group $A$, define the \textbf{$p$-primary subgroup of $A$} to be 
$$\empty_{p^\infty}(A)=\{a \in A ~|~ a^{p^i}=1 \text{ for some }i >0 \}.$$
The order of this subgroup  is of the form $p^e$. Call $e$ the
\textbf{$p^\infty$-rank of $A$}.
\end{definition} 
\noindent The $p^\infty$ rank of a finitely-generated 
abelian group $A$ can be calculated by taking as input the abelian invariants of $A$ and the 
prime $p$.

By passing to $p$-primary subgroups, sequence \ref{8} gives another
short exact sequence
\[
1\to k\otimes\frac{R[F,F]}{[F,F]}\to
\empty_{p^\infty} \left( \frac{F}{R^p[F,F]} \right) \to
\empty_{p^\infty} \left( \frac{F}{R[F,F]} \right) \to 1
\]
since the first term is $p$-torsion.
We observe that while $F/R[F,F]$ can be given in terms of $S(F)$ and
$S(R)$, the factor group $F/R^p[F,F]$ can be given in the same way
but replacing $S(R)$ by $S(R)^p$ - the finite list of $p$-powers of
elements in $S(R)$.\\

\subsection{The Rank Algorithm}
~\\

\noindent \textbf{Algorithm 4:} $\proc{PrimePrimaryRank}(F, R, p)$
\begin{algorithmic}[1]
\REQUIRE Free Group $F$, Relators $R$, Prime $p$
\ENSURE $p^\infty$-rank of $F/R$
\STATE $A:=$FirstHomology$(F,R)$
\STATE $Y:=[\ ]$
\FOR{$a \in A$}
\IF{$a \ne 0$ and $a \equiv 0 \text{ mod }p$}
\STATE $y:=$ p-adiv valuation of $a$
\STATE append $y$ to $Y$
\ENDIF
\ENDFOR
\STATE $s:=$Sum$(Y)$ \COMMENT{$s$ is the sum of the elements of $Y$}
\RETURN $s$
\end{algorithmic}
~\\

The GAP command \verb+PadicValuation(n,p)+ gives the $p$-adic valuation of an integer $n$.

To summarize, let
\begin{eqnarray*}
a&=& \text{ dimension of Tor}(H_1(G),k)\\
b&=& p^\infty\text{-rank of }\frac{F}{R[F,F]}\\
c&=& p^\infty\text{-rank of }\frac{F}{R^p[F,F]}\\
d&=& \text{ dimension of }H_2(G;k)\\
e&=& \text{ dimension of }k\otimes\frac{R}{[F,R]}
\end{eqnarray*}
where $a$ is determined by the Tor Algorithm, $b$ and $c$ by the
Rank Algorithm, and $e$ is yet to be studied. By the additive
property of the dimension and the $p^\infty$-rank we deduce, from the
exact sequences above, the following reduction formula:
\[
d=a+b-c+e.
\]
Since $a$, $b$, $c$ are more or less standard, the integer $e$ is the
key difficulty we aim to approach experimentally.

We first describe an algorithm that reduces an element of a group via a rewriting system.\\

\subsection{Reduce Word Algorithm}
~\\

\noindent \textbf{Algorithm 5:} $\proc{Reduce\_Word}(F, R, Z, R', p)$
\begin{algorithmic}[1]
\REQUIRE Free Group $F$, Relators $R$, Test Word $z$, Sublist $R'$ of $R$, Prime $p$
\ENSURE Reduced word of $z$ in $F/[F,R]R^pR'$
\STATE $G:=F/[F,R]R^pR'$
\STATE $RG:=$Rewriting system for $G$
\STATE $x:=$ReducedWord$(z)$
\RETURN $x$
\end{algorithmic}
~\\

We use the rewriting system given by the Knuth-Bendix completion algorithm \cite{knuth:70} implemented
on GAP via the KBMAG package \cite{holt-kbmag:www}.\\

\subsection{The Find Basis Algorithm}
~\\

\noindent \textbf{Algorithm 6:} $\proc{FindBasis}(F,R,p,R')$
\begin{algorithmic}[1]
\REQUIRE Free Group $F$, Relators $R$, Prime $p$, Sublist $R'$ of $R$
\ENSURE Size of a generating set for $[F,R]R^pR'/[F,R]R^p$
\STATE $X:=R'$
\FOR{$x \in X$}
\STATE $x':=$Reduce$\_$Word$(F,R,x,\text{Difference}(X,[x]),p)$ \COMMENT{Difference$(A,B)$ is the complement of $B$ in $A$}
\IF{$x'=$ identity}
\STATE $X:=$Difference$(X,[x])$
\ENDIF
\ENDFOR
\RETURN Size$(X)$
\end{algorithmic}
~\\

The algorithm attempts to check for linear independence of each element $x$ of $R'$ with 
respect to $R'-\{x\}$ in $[F,R]R^pR'/[F,R]R^p$. Whenever $x$ is found by the rewriting system
to be dependent of $R'-\{x\}$, it is removed from $R'$. The end result will be a list of potentially 
linearly independent generators.

We conclude this discussion with the grand scheme algorithm
which takes as input a finite list of generators $S(F)$ and a finite
list of relators $S(R)$ for a group $G$ together with a prime $p$
and gives as output an integer $d$ representing an upper bound for
the dimension of $H_2(G;k)$, where $k$ is a field of characteristic $p$.\\

\subsection{The Second Homology with Coefficients Algorithm}
~\\

\noindent \textbf{Algorithm 7:} $\proc{SecondHomologyCoefficients}(F, R, p, R')$
\begin{algorithmic}[1]
\REQUIRE Free Group $F$, Relators $R$, Prime $p$, Sublist $R'$ of $R$ generating $R/[F,R]R^p$
\ENSURE An integer $d$ such that dim $\left(H_2(G;k)\right) \le d$
\STATE $a:=$ Tor$(F,R,p)$
\STATE $b:=$ PrimePrimaryRank$(F,R[F,F],p)$
\STATE $c:=$ PrimePrimaryRank$(F,R^p[F,F],p)$
\STATE $e:=$ FindBasis$(F,R,p,R')$
\STATE $d:=a+b-c+e$
\RETURN $d$
\end{algorithmic}
~\\

It is important to note that the reduction of test words in the algorithm
``Reduce\_ Word'' is the word problem. As such, a result of a word not being the 
identity is an indeterminate result. However, if $G$ is finite, or, more generally, 
if the rewriting is confluent, the reduction in the rewriting system 
is deterministic and a basis is achieved (the confluence for finite groups is guaranteed in theory only; in practice
it may take a long time or require more space than is available). At any rate, this is not typically 
the case--the word problem is undecidable in general, thus the result 
of ``Find\_ Basis'' is, in general, the cardinality of a generating set that is not necessarily a basis. Therefore
in these cases we do not find the dimension of $H_2(G;k)$, only an upper bound.

\section{Examples}\label{examples}
In this section, we apply the grand scheme algorithm above to some select groups. The first example is to 
illustrate the effect the algorithm has on groups with smallish presentations. The other three examples 
are the groups of primary interest.

\begin{example} \label{ex1} The symmetric groups $\Sigma_n$ on $n$ letters \cite{brown:82}.
\begin{eqnarray*}
G&=&\Sigma_5\\
S(F)&=&[a,b]\\
S(R)&=&[a^5,b^2,(a^{-1}b)^4,(a^2ba^{-2}b)^2]\\
p&=&2\\
d&=&2
\end{eqnarray*}
\end{example}

Next, we consider three linear groups over $\mathbb Z[1/p,\zeta_p]$ where $\zeta_p$ is a
primitive $p^\text{th}$-root of unity. Presentations for groups of this from can be found in \cite{anton:08}.
\begin{example} \label{ex2} \item \begin{eqnarray*}
G&=&SL_2(\mathbb Z[1/3,\zeta_3])\\
S(F)&=&[z, u_1, a, b, b_0, b_1, b_2, w]\\
S(R) &=&[b_t^{-1}z^{3t}bz^{3t}a,w^{-1}z^4u_1u_2u_3,z^3, [z,u_1], [u_1,u_1],a^4, [a^2,z], [a^2,u_1],\\
&& a^{-1}zaz, a^{-1}u_1au_1,\left[ b_s, b_t \right], b^{-3}a^2, b^{-3}b_0b_1b_2,\\
&& (b_0b_1^{-1}a^{-1}u_1)^3, a^{-2}b^{-1}u_1bz^{-3}b^{-1}b_0^{-1}z^{3}bz^{-1}u_1]\\
p&=& 3\\
d&=& 0
\end{eqnarray*}
where $s,t \in \{1,2\}$.\\
\end{example}

\begin{example} \label{ex3}
\begin{eqnarray*}
G&=&SL_2(\mathbb Z[1/5,\zeta_5])\\
S(F)&=&[z, u_1, u_2, a, b, b_0, b_1, b_2, b_3, b_4, w]\\
S(R)&=&[b_t^{-1}z^{3t}bz^{3t}a,w^{-1}z^4u_1u_2u_3,z^5, [z,u_i], [u_i,u_j],a^4, [a^2,z], [a^2,u_i],\\
&& a^{-1}zaz, a^{-1}u_iau_i,\left[ b_s, b_t \right], b^{-3}a^2, b^{-3}b_0b_1b_2b_3b_4,\\
&& (b_0b_1^{-1}a^{-1}u_1)^3, (b_0b_2^{-1}a^{-1}u_2)^3, (b_0b_3^{-1}a^{-1}u_3)^3,\\
&& (b_0b_1^{-1}b_2^{-1}b_3a^{-1}u_1u_2)^3, (b_0b_1^{-1}b_3^{-1}b_4a^{-1}u_1u_3)^3,\\
&& (b_0b_2^{-1}b_3^{-1}b_5a^{-1}u_2u_3)^3,  a^{-2}b^{-1}u_ibz^{-3i}b^{-1}b_0^{-1}z^{3i}bz^{-i}u_i]\\
p&=&5\\
d&=&0
\end{eqnarray*}
where $i,j \in \{1,2\}$ and $s,t \in \{1,2,3,4\}$.
\end{example}

\begin{example} \label{ex4}
\begin{eqnarray*}
G&=&SL_2(\mathbb Z[1/7,\zeta_7])\\
S(F)&=&[z,u_1,u_2,u_3,a,b,b_0,b_1,b_2,b_3,b_4,b_5,b_6,w]\\
S(R)&=&[b_t^{-1}z^{3t}bz^{3t}a,w^{-1}z^4u_1u_2u_3,z^7, [z,u_i], [u_i,u_j],a^4, [a^2,z], [a^2,u_i],\\
&&a^{-1}zaz, a^{-1}u_iau_i,\left[ b_s, b_t \right], b^{-3}a^2, b^{-3}b_0b_1b_2b_3b_4b_5b_6,b^{-7}_tw^{-1}b_t^{-1}w,\\
&&(b_0b_1^{-1}a^{-1}u_1)^3,(b_0b_2^{-1}a^{-1}u_2)^3,(b_0b_3^{-1}a^{-1}u_3)^3,\\
&&(b_0b_1^{-1}b_2^{-1}b_3a^{-1}u_1u_2)^3, (b_0b_1^{-1}b_3^{-1}b_4a^{-1}u_1u_3)^3,(b_0b_2^{-1}b_3^{-1}b_5a^{-1}u_2u_3)^3,\\
&&(b_0b_1^{-1}b_2^{-1}b_3b_4b_5b_6^{-1}a^{-1}u_1u_2u_3)^3, a^{-2}b^{-1}u_ibz^{-3i}b^{-1}b_0^{-1}z^{3i}bz^{-i}u_i]\\
p&=&7\\
d&= &6
\end{eqnarray*}
where $i,j \in \{1,2,3\}$ and $s,t \in \{1,2,3,4,5,6\}$.
\end{example}

\section{Discussion and Future Work}

\noindent Details on the above examples are as follows:\\

\begin{itemize}
\item Example \ref{ex1}: The rewriting system given by the KBMAG package for $\Sigma_5$ is confluent; therefore
 $$\text{dim} H_2(\Sigma_5;\mathbb F_2)=2.$$
 The algorithm took about 50 milliseconds to run,
reflecting the relatively simple presentation. 

\item Example \ref{ex2}: The rewriting system given by the KBMAG package for $SL_2(\mathbb Z[1/3,\zeta_3])$ is not confluent; the algorithm
took about 6 hours to finish. In this case, the non-confluence of the system did not affect the 
results as the rewriting system was able to show that all elements of $R$ reduced to identity modulo $[F,R]R^3$, 
so $$\text{dim} H_2(SL_2(\mathbb Z[1/3, \zeta_3]; \mathbb F_3)=0.$$
 
\item Example \ref{ex3}: The rewriting system given by the KBMAG package for $SL_2(\mathbb Z[1/5,\zeta_5])$ is not confluent. As in Example 2
the non-confluence of the system did not affect the results and 
$$\text{dim} H_2(SL_2(\mathbb Z[1/5, \zeta_5])=0.$$
The algorithm took about 2 days to finish.

\item Example \ref{ex4}: The rewriting system given by the KBMAG package for $SL_2(\mathbb Z[1/7,\zeta_7])$ is not confluent. In this case,
the algorithm took a total of about 5 days to finish. Also, this is the only case tested in which the 
non-confluence actually mattered. Since the algorithms were not able to show that the dimension of $R/[F,R]R^7$ is 0,
we only have the upper bound 
$$\text{dim} H_2(SL_2(\mathbb Z[1/7, \zeta_7]); \mathbb F_7) \le 6.$$
\end{itemize}

We note that for Examples \ref{ex3} and \ref{ex4}, it was necessary to run the algorithm several times to
obtain the results above since the parameters of the KBMAG package allow a limited number
of equations to be generated in the rewriting system. Each iteration eliminated elements of $R$ from the 
generating list until the results stabilized. For instance, in Example \ref{ex4}, the initial iteration
gave a result of $e \le 16$ and $d \le 10$, the second iteration gave that $e \le 13$ and $d \le 7$.
The third and fourth iterations each gave a result of $e \le 12$ and so the upper bound on $d$ is $6$.

Finally, in implementing these algorithms to find a bound on $H_2(G)$ it is useful
to first perform Tietze transforms on the presentations involved to attempt to simplify
the presentations. In many cases, the number of generators and relators can be 
reduced, thus simplifying the calculations. In Example \ref{ex4} 
$SL_2(\mathbb Z[1/7, \zeta_7])$ is given via a presentation consisting of 14
generators and 64 relators. A series of Tietze transforms, implemented via GAP,
simplifies to a presentation with 6 generators and 34 relators. This significantly 
impacts the results of the algorithm.

Our future work will involve refining and improving the algorithms above. Initially we 
were concerned only with writing algorithms that gave results--the efficiency 
of these algorithms was not a concern. For the linear groups above as $p$ increases the number of relators grows exponentially, thus the
algorithms will take longer and longer to finish. Going from $p=2$ to $p=7$ the time required
increased from several hours to several days.

We also will develop other methods for finding generators of $H_2(G)$ and $H_2(G;k)$ independent
from those above. In particular, we attempt to find \textit{lower} bounds on the dimension of $H_2(G;k)$.
The strategies for both problems will be based on linear algebra involving rewriting systems for $S(R)$ 
in $F/[F,F]$ and will appear in a future paper.

\appendix
\section{Appendix of functions} \label{appendix}

In this appendix we give the code used in GAP to implement the algorithms described above.

\begin{verbatim}
##############################################
#Input: Free group, relators, prime p
#Output: Dimension of Tor(H_1(G),F_p)

Tor:=function(Freegroup,Relators,Prime)
local AbelInv, list1, list2, x;
AbelInv:=AbelianInvariants(Freegroup/Relators);
list1:=[];
list2:=[];
for x in AbelInv do
     if x<>0 then Add(list1,x)
     fi;
     od;
for x in list1 do
     if x mod Prime = 0 then Add(list2,1)
     fi;
     od;
return Sum(list2);
end;;

###############################################
#P-Primary Rank
#Output: P-Primary Rank of Fp Group

PrimePrimaryRank:=function(Freegroup,Relators,Prime)
local AbelInv, list1, list2, x;
AbelInv:=AbelianInvariants(Freegroup/Relators);
list1:=[];
list2:=[];
for x in AbelInv do
     if x <> 0 then Add(list1,x)
     fi;
     od;
for x in list1 do
     if x mod Prime = 0 then Add(list2,PadicValuation(x,Prime))
     fi;
	     od;
return Sum(list2);	
end;;

###############################################
#The (special) Word Problem
#Output: Reduced word

Reduce_Word:=function(Freegroup,Relators,TestWord,Sublist,Prime)
local Rel_P, GroupGen, comm, G, RG, OR;

Rel_P:=List(Relators,x->x^Prime);
GroupGen:=GeneratorsOfGroup(Freegroup);
comm:=ListX(GroupGen,Relators,Comm);

G:=Freegroup/Concatenation(comm,Rel_P,Sublist);
RG:=KBMAGRewritingSystem(G);
OR:=OptionsRecordOfKBMAGRewritingSystem(RG);
OR.maxeqns:=500000;
OR.tindyint:=100;
MakeConfluent(RG);

return ReducedWord(RG,TestWord);
end;;

#################################################
#Attempts to reduce a generating set
#Output: List of generators

FindBasis:=function(Freegroup,Relators,Prime,Sublist)
local Gen,TestWord,x; 
Gen:=Sublist;
for x in Sublist do 
     TestWord:=Reduce_Word(Freegroup,Relators,x,Difference(Gen,[x]),Prime)
        if IsOne(TestWord)=true then Gen:=Difference(Gen,[x])
        fi;
        od;
return Size(Gen);
end;;

#######################################################
#Gives the estimate for H_2
#Output: Upper bound on dimension of H_2
SecondHomologyCoefficients:=function(Freegroup, Relators, Prime, Sublist)

local a,b,c,d,e,f,ff,RPrime;

f:=GeneratorsOfGroup(Freegroup);
ff:=ListX(f,f,Comm);
RPrime:=List(Relators,x->x^Prime);

a:=Tor(Freegroup,Relators,Prime);
b:=PrimePrimaryRank(Freegroup,Concatenation(Relators,ff),Prime);
c:=PrimePrimaryRank(Freegroup,Concatenation(RPrime,ff),Prime);
e:=FindBasis(Freegroup,Relators,Prime,Sublist);
d:=a+b-c+e;
return d;
end;;
\end{verbatim}

\bibliographystyle{amsplain}
\bibliography{ref}

\end{document}